\documentclass[11pt]{article}
\usepackage{amscd}
\usepackage{amsfonts}
\usepackage{amsmath}
\usepackage{amssymb}
\usepackage{amsthm}
\usepackage{bbm}
\usepackage{CJK}
\usepackage{fancyhdr}
\usepackage{graphicx}
\usepackage{hyperref}
\usepackage{indentfirst}
\usepackage{latexsym}
\usepackage{mathrsfs}
\usepackage{xypic}
\usepackage{amsmath,amssymb,amsfonts,amsthm,fancyhdr}
\usepackage{epsfig,graphicx,picins,picinpar,subfigure}
\usepackage{pstricks}
\usepackage{amssymb}
\usepackage{xypic}
\usepackage[compress]{cite}

\newtheorem{theorem}{Theorem}[section]
\newtheorem{prop}[theorem]{Proposition}
\newtheorem{lemma}[theorem]{Lemma}

\newtheorem{prop-def}{Proposition-Definition}[section]

\theoremstyle{definition}
\newtheorem{defn}[theorem]{Definition}

\newtheorem{remark}[theorem]{Remark}
\newtheorem{exam}[theorem]{Example}

\usepackage[top=1in,bottom=1in,left=1.25in,right=1.25in]{geometry}
\def\<{\langle}
\def\>{\rangle}
\def\a{\alpha}

\date{\today}
\begin{document}
\renewcommand{\baselinestretch}{1.2}
\renewcommand{\arraystretch}{1.0}
\title{\bf  On  compatible   Hom-Lie triple systems }
\author{{\bf Wen Teng$^{1}$,~   Fengshan Long$^1$,~ Hui Zhang$^{2,3}$,~ Jiulin Jin$^{4}$}\\
{\small 1. School of Mathematics and Statistics, Guizhou University of Finance and Economics,} \\
{\small  Guiyang  550025, P. R. of China}\\
  {\small Email: tengwen@mail.gufe.edu.cn (Wen Teng)} \\
{\small 2.   School of Information, Guizhou University of Finance and Economics,   }\\
{\small Guiyang  $550005$, P. R. of China}\\
{\small 3. Postdoctoral Scientific Research Station, ShijiHengtong Technology Co.,} \\
{\small   Ltd, Guiyang,  550014, P. R. of China}\\
{\small 4.   School of Science,  Guiyang University,   Guiyang  550005,  China}
}

 \maketitle
\begin{center}
\begin{minipage}{13.cm}

{\bf Abstract}
In this paper, we consider compatible Hom-Lie triple systems.
Compatible Hom-Lie triple systems are characterized as Maurer-Cartan elements in a suitable bidifferential
graded Lie algebra. We also define a cohomology theory for compatible Hom-Lie triple systems. As applications of cohomology, we study
deformations and  abelian extensions of compatible  Hom-Lie triple systems.

 \smallskip

 \smallskip

 \end{minipage}
 \end{center}
 \normalsize\vskip0.5cm

\section{Introduction}
\def\theequation{\arabic{section}. \arabic{equation}}
\setcounter{equation} {0}

Lie triple system  first appeared in Cartan's work \cite{Cartan} on Riemannian geometry. Since then,
In \cite{Jacobson1,Jacobson2}, Jacobson studied Lie triple systems from Jordan theory and quantum mechanics. After that,
Lie triple system has important applications in physics, such as quantum mechanics theory, elementary particle theory and numerical analysis of differential equations.
The representation  theory, deformation theory, cohomology and related properties of Lie triple system can be found in references \cite{Yamaguti,  Kubo}.
As a Hom-type  algebra \cite{Hartwig} generalization of Lie triple system, Hom-Lie triple system was introduced by Yau in \cite{Yau}.
Later,   Ma et al. \cite{Ma} considered   cohomology and 1-parameter formal deformations of Hom-Lie triple systems. The authors have studied the relative Rota-Baxter operators on Hom-Lie triple systems in \cite{Teng3, Teng4,Li}.

The purpose of this paper is to study the compatible Hom-Lie triple system, which is a pair of Hom-Lie triple systems, so that the linear combination of their algebraic structures is also Hom-Lie triple system.
Compatible algebraic structures are widely used in various fields of mathematics and mathematical physics.
For example, in the study of bi-Hamiltonian mechanics, the concept of compatibility of two Poisson structures on manifolds is introduced   from the mathematical point of view, see \cite{Magri,Kosmann}.
Recently, compatible algebraic structures have been widely studied, such as  compatible Lie algebra \cite{Liu23},  compatible $L_\infty$-algebras \cite{Das22}, compatible pre-Lie algebra \cite{Liu}, compatible associative algebras \cite{Chtioui},     compatible Lie-Yamaguti algebras \cite{Sania}, compatible Lie supertriple system \cite{Wang}, compatible 3-Lie algebras \cite{Hou},  compatible Hom-Lie algebras \cite{Das23} and compatible Hom-Leibniz algebras \cite{Bhutia}.
Motivated by these works, in this paper, we study the cohomology,   deformations and abelian extensions of compatible Hom-Lie triple systems.

The paper is organized as follows.   In Section  \ref{sec: compatible Hom-Lie triple system},  we introduce the notion of compatible Hom-Lie triple systems and give the
bidifferential graded Lie algebra whose Maurer-Cartan elements are compatible Hom-Lie triple system structures.
In Section \ref{sec: Cohomologies},   we introduce the cohomology of a compatible Hom-Lie triple system with coefficients in a representation.
 In Section \ref{sec: deformations},  we study linear deformations  of a compatible Hom-Lie triple system. We introduce Nijenhuis operators that
generate trivial linear deformations.
 In Section \ref{sec: abelian extensions},  we study abelian extensions of compatible
Hom-Lie triple systems and give a classification of equivalence classes of abelian extensions.

Throughout this paper, $\mathbb{K}$ denotes a field of characteristic zero. All the    vector spaces  and
   (multi)linear maps are taken over $\mathbb{K}$.

\section{ Maurer-Cartan characterizations of compatible Hom-Lie triple systems }\label{sec: compatible Hom-Lie triple system}
\def\theequation{\arabic{section}.\arabic{equation}}
\setcounter{equation} {0}

In this section,  we  recall  concepts of    Hom-Lie triple systems from \cite{Yau} and \cite{Ma}.
  Then we introduce  the notion of compatible Hom-Lie triple systems and give the
bidifferential graded Lie algebra whose Maurer-Cartan elements are compatible Hom-Lie triple system
structures. We also construct the bidifferential graded Lie algebra governing deformations of a compatible
Hom-Lie triple system.

\begin{defn}  \cite{Yau}
(i) A Hom-Lie triple system (Hom-Lts)  is a vector space $\mathfrak{g}$ together with  a trilinear operation $[-, -, -]_\mathfrak{g}$
on $\mathfrak{g}$  and a linear map $\alpha_\mathfrak{g}: \mathfrak{g} \rightarrow  \mathfrak{g}$,    satisfying $\alpha_\mathfrak{g}([x,y,z]_\mathfrak{g})=[\alpha_\mathfrak{g}(x),\alpha_\mathfrak{g}(y),\alpha_\mathfrak{g}(z)]_\mathfrak{g}$ such that
\begin{align}
&[x,y,z]_\mathfrak{g}+[y,x,z]_\mathfrak{g}=0,\label{2.1}\\
&[x,y,z]_\mathfrak{g}+[z,x,y]_\mathfrak{g}+[y,z,x]_\mathfrak{g}=0,\label{2.2}\\
 &[\alpha_\mathfrak{g}(a), \alpha_\mathfrak{g}(b), [x, y, z]_\mathfrak{g}]_\mathfrak{g}=[[a, b, x]_\mathfrak{g}, \alpha_\mathfrak{g}(y), \alpha_\mathfrak{g}(z)]_\mathfrak{g}+ [\alpha_\mathfrak{g}(x),  [a, b, y]_\mathfrak{g}, \alpha_\mathfrak{g}(z)]_\mathfrak{g}\nonumber\\
 &~~~~~~~~~~~~~~~~~~~~~~~~~~~~~~~~~~~~~~~~~~~~~~~~~~~~~~~~~~~~~~~+ [\alpha_\mathfrak{g}(x),  \alpha_\mathfrak{g}(y), [a, b, z]_\mathfrak{g}]_\mathfrak{g},\label{2.3}
\end{align}
where $ x, y, z, a, b\in \mathfrak{g}$.  A Hom-Lie triple system  as above  may be denoted by the triple  $(\mathfrak{g}, [-, -, -]_\mathfrak{g},\alpha_\mathfrak{g})$.
In particular, the Hom-Lie triple system  $(\mathfrak{g}, [-, -, -]_\mathfrak{g},\alpha_\mathfrak{g})$  is said to be regular, if $\alpha_\mathfrak{g}$ is nondegenerate.

(ii) A homomorphism between two  Hom-Lie triple systems  $(\mathfrak{g}_1, [-, -, -]_1,\alpha_1)$ and $(\mathfrak{g}_2, [-, -,$ $-]_2,\alpha_2)$ is a linear map $\varphi: \mathfrak{g}_1\rightarrow \mathfrak{g}_2$ satisfying
$$\varphi(\alpha_1(x))=\alpha_2(\varphi(x)), ~~\varphi([x, y, z]_1)=[\varphi(x), \varphi(y),\varphi(z)]_2,~~\forall ~x,y,z\in \mathfrak{g}_1.$$
\end{defn}

\begin{remark}
(i) Let $(\mathfrak{g}, [-,-]_\mathfrak{g},  \alpha_\mathfrak{g})$ be a     Hom-Lie   algebra, then
$( \mathfrak{g}, [-,-,-]_\mathfrak{g},  \alpha_\mathfrak{g})$ is a  Hom-Lie triple system,
where  $ [x,y,z]_{\mathfrak{g}}=[[x,y]_\mathfrak{g},\alpha_\mathfrak{g}(z)]_\mathfrak{g}, \forall x,y,z\in \mathfrak{g}.$

(ii)  A  Hom-Lie triple system   $(\mathfrak{g}, [-,-]_\mathfrak{g},  \alpha_\mathfrak{g})$  with   $\alpha_\mathfrak{g}=\mathrm{id}_\mathfrak{g}$  is nothing but a Lie triple system,
 it follows that the results established in this paper can be naturally adapted to the context
of Lie triple system.
\end{remark}

\begin{defn} \cite{Ma}
A representation of  a  Hom-Lie triple system  $(\mathfrak{g}, [-, -,-]_\mathfrak{g},\alpha_\mathfrak{g})$ on  a Hom-vector space $(V,\beta)$ is a bilinear map $\theta: \mathfrak{g}\times \mathfrak{g}\rightarrow \mathrm{End}(V)$, such that for all  $x,y,a,b\in \mathfrak{g}$
\begin{align}
& \theta(\alpha_\mathfrak{g}(x),\alpha_\mathfrak{g}(y))\circ\beta=\beta\circ\theta(x,y),\label{2.4}\\
& \theta(\alpha_\mathfrak{g}(a),\alpha_\mathfrak{g}(b))\theta(x,y)-\theta(\alpha_\mathfrak{g}(y),\alpha_\mathfrak{g}(b))\theta(x,a)-\theta(\alpha_\mathfrak{g}(x), [y,a,b]_\mathfrak{g})\circ\beta\nonumber\\
&+D(\alpha_\mathfrak{g}(y),\alpha_\mathfrak{g}(a))\theta(x,b)=0,\label{2.5}\\
& \theta(\alpha_\mathfrak{g}(a),\alpha_\mathfrak{g}(b))D(x,y)-D(\alpha_\mathfrak{g}(x),\alpha_\mathfrak{g}(y))\theta(a,b)+\theta([x,y,a]_\mathfrak{g},\alpha_\mathfrak{g}(b))\circ\beta\nonumber\\
&+\theta(\alpha_\mathfrak{g}(a),[x,y,b]_\mathfrak{g})\circ\beta=0,\label{2.6}\
\end{align}
 where $D(x,y)=\theta(y,x)-\theta(x,y)$. We also denote a representation of $\mathfrak{g}$ on $(V,\beta)$ by $(V,\beta; \theta)$.
\end{defn}

\begin{prop} \cite{Ma}
Let $(V,\theta; \beta)$   be a representation of a  Hom-Lie triple system  $(\mathfrak{g},[-,-,-]_\mathfrak{g}, \alpha_\mathfrak{g})$.
Define a trilinear operation $[-,-,-]_{\ltimes}:\wedge^2(\mathfrak{g}\oplus V)\otimes(\mathfrak{g}\oplus V)\rightarrow  \mathfrak{g}\oplus V$
and a linear map $\alpha_\mathfrak{g}\oplus\beta: \mathfrak{g}\oplus V \rightarrow  \mathfrak{g}\oplus V$ by
\begin{align*}
[(a,u),(b,v),(c,w)]_{\ltimes}&=([a,b,c]_\mathfrak{g},D(a,b)w+\theta(b,c)u-\theta(a,c)v),\\
\alpha_\mathfrak{g}\oplus\beta(a,u)&=(\alpha_\mathfrak{g}(a),\beta(u)),
\end{align*}
for any $(a,u),(b,v),(c,w)\in \mathfrak{g}\oplus V$.  Then $(\mathfrak{g}\oplus V,[-,-,-]_{\ltimes}, \alpha_\mathfrak{g}\oplus\beta)$ is a  Hom-Lie triple system.
This  Hom-Lie triple system is called the semi-product Hom-Lie triple system and denoted by  $\mathfrak{g}\ltimes V$.
\end{prop}

Next we recall the cohomology theory on Hom-Lie triple systems given in  \cite{Ma}.
Let $(V,\beta; \theta)$  be a representation of a  Hom-Lie triple system $( \mathfrak{g}, [-,-,-]_\mathfrak{g},\alpha_\mathfrak{g})$.
Denote the $n$-cochains   of $ \mathfrak{g}$ with coefficients in representation $(V,\beta; \theta)$   by
\begin{align*}
%\mathcal{C}_{\mathrm{HLts}}^{0}( \mathfrak{g},  V):&=\{a\in \mathfrak{g}~|~\alpha_\mathfrak{g}(a)=a\}, \\
\mathcal{C}_{\mathrm{HLts}}^{1}( \mathfrak{g},  V):&=\{f\in  \mathrm{Hom}(\mathfrak{g},V)~|~\beta(f(a))=f(\alpha_\mathfrak{g}(a))\},\\
\mathcal{C}_{\mathrm{HLts}}^{n+1}( \mathfrak{g},  V):&=\big\{f\in  \mathrm{Hom}  (\overbrace{\wedge^2 \mathfrak{g}\otimes\cdots\otimes\wedge^2 \mathfrak{g}}^n\otimes  \mathfrak{g},V) ~|~ \beta\circ f=f\circ (\alpha_\mathfrak{g}^{\wedge2}\otimes\cdots\otimes\alpha_\mathfrak{g}^{\wedge2}\otimes\alpha_\mathfrak{g}),\\
&f(a_1\wedge b_1,\cdots, a_{n-1}\wedge b_{n-1},a,b,c)+f(a_1\wedge b_1,\cdots, a_{n-1}\wedge b_{n-1},b,a,c)=0,\\
&\circlearrowright_{a,b,c}f(a_1\wedge b_1,\cdots, a_{n-1}\wedge b_{n-1},a,b,c)=0  \big\}.
\end{align*}

The   coboundary operator $\delta: \mathcal{C}_{\mathrm{HLts}}^{n}(\mathfrak{g}, V)\rightarrow \mathcal{C}_{\mathrm{HLts}}^{n+1}( \mathfrak{g}, V)$ by
\begin{align*}
&\delta f( \mathfrak{A}_1,\cdots,  \mathfrak{A}_{n},c)\\
=&\theta(\alpha^{n-1}_{\mathfrak{g}}(b_{n}),\alpha^{n-1}_{\mathfrak{g}}(c))f(\mathfrak{A}_1,\cdots,  \mathfrak{A}_{n-1} a_{n})-\theta(\alpha^{n-1}_{\mathfrak{g}}(a_{n}),\alpha^{n-1}_{\mathfrak{g}}(c))f(\mathfrak{A}_1,\cdots,  \mathfrak{A}_{n-1},b_{n})\\
&+\sum_{i=1}^n(-1)^{i+n}D(\alpha^{n-1}_{\mathfrak{g}}(a_{i}),\alpha^{n-1}_{\mathfrak{g}}(b_{i}))f(\mathfrak{A}_1,\cdots,  \widehat{\mathfrak{A}}_{i},\cdots,  \mathfrak{A}_{n}\cdots, c)\\
&+\sum_{i=1}^n\sum_{j=i+1}^{n}(-1)^{i+n+1}f(\alpha_{\mathfrak{g}}( \mathfrak{A}_1),\cdots,   \widehat{\mathfrak{A}}_{i}, \cdots, [a_{i},b_{i},a_{j}]_\mathfrak{g}\wedge \alpha_{\mathfrak{g}}(b_j)+\alpha_{\mathfrak{g}}(a_j)\wedge [a_{i},b_{i},b_{j}]_\mathfrak{g},\cdots,\\
&~~~~~~~~~ \alpha_{\mathfrak{g}}(\mathfrak{A}_{n}),\alpha_{\mathfrak{g}}(c))+\sum_{i=1}^n(-1)^{i+n+1}f(\alpha_{\mathfrak{g}}( \mathfrak{A}_1),\cdots,  \widehat{\mathfrak{A}}_{i},\cdots,     \alpha_{\mathfrak{g}}(\mathfrak{A}_{n}), [a_{i},b_{i},c]_\mathfrak{g}),
\end{align*}
where   $  \mathfrak{A}_i=a_i\wedge b_{i}\in  \wedge^2 \mathfrak{g}, 1\leq i\leq n, c\in \mathfrak{g} $,  $f\in \mathcal{C}_{\mathrm{HLts}}^{n}( \mathfrak{g}, V)$,  $\alpha_{\mathfrak{g}}( \mathfrak{A}_i)=\alpha_{\mathfrak{g}}(a_i)\wedge\alpha_{\mathfrak{g}}(b_i).$

 Let $(\mathfrak{g},\alpha_\mathfrak{g})$ be a Hom-vector space. For each $n\geq 1$, consider the spaces  $C^n_\mathrm{H}(\mathfrak{g},\mathfrak{g})$  by
\begin{align*}
%C_{\mathrm{H}}^{0}( \mathfrak{g},  \mathfrak{g}):&=\{a\in \mathfrak{g}~|~\alpha_\mathfrak{g}(a)=a\}, \\
C_{\mathrm{H}}^{1}( \mathfrak{g},  \mathfrak{g}):&=\{f\in  \mathrm{Hom}(\mathfrak{g},\mathfrak{g})~|~\alpha_\mathfrak{g}(f(a))=f(\alpha_\mathfrak{g}(a))\},\\
C_{\mathrm{H}}^{n+1}( \mathfrak{g}, \mathfrak{g}):&=\big\{f\in  \mathrm{Hom}  (\overbrace{\wedge^2 \mathfrak{g}\otimes\cdots\otimes\wedge^2 \mathfrak{g}}^n\otimes  \mathfrak{g},\mathfrak{g}) ~|~\alpha_\mathfrak{g}\circ f=f\circ (\alpha_\mathfrak{g}^{\wedge2}\otimes\cdots\otimes\alpha_\mathfrak{g}^{\wedge2}\otimes\alpha_\mathfrak{g}),\\
&f(a_1\wedge b_1,\cdots, a_{n-1}\wedge b_{n-1},a,b,c)+f(a_1\wedge b_1,\cdots, a_{n-1}\wedge b_{n-1},b,a,c)=0 \big\}.
\end{align*}

Then the   graded vector space $C^{*+1}_\mathrm{H}(\mathfrak{g},\mathfrak{g})=\oplus_{n=0}^{+\infty}C^{n+1}_\mathrm{H}(\mathfrak{g},\mathfrak{g})$  carries a graded Lie bracket defined as
follows. For $P\in C^{p+1}_\mathrm{H}(\mathfrak{g},\mathfrak{g})$ and  $Q\in C^{q+1}_\mathrm{H}(\mathfrak{g},\mathfrak{g})$, the   bracket  $[P,Q]_{\mathrm{Hlts}}\in C^{p+q+1}_\mathrm{H}(\mathfrak{g},\mathfrak{g})$
given by
\begin{align*}
 &[P,Q]_{\mathrm{Hlts}}=P\diamond Q-(-1)^{pq}Q\diamond P, ~~  ~~\text{where} \\
 &P\diamond Q(\mathfrak{A}_1,\cdots,\mathfrak{A}_{p+q},c)\\
= &\sum_{\sigma\in \mathbb{S}(p,q)}(-1)^{pq+\sigma}P(\alpha^n_{\mathfrak{g}}(\mathfrak{A}_{\sigma(1)}),\cdots,\alpha^n_{\mathfrak{g}}(\mathfrak{A}_{\sigma(p)}), Q(\mathfrak{A}_{\sigma(p+1)},\cdots,\mathfrak{A}_{\sigma(p+q)},c))+\\
 &\sum_{k=1}^p\sum_{\sigma\in \mathbb{S}(k-1,q)}(-1)^{(k-1)q+\sigma}P(\alpha^n_{\mathfrak{g}}(\mathfrak{A}_{\sigma(1)}),\cdots,\alpha^n_{\mathfrak{g}}(\mathfrak{A}_{\sigma(k-1)}), \alpha^n_{\mathfrak{g}}(a_{k+q})\wedge Q(\mathfrak{A}_{\sigma(k)},\cdots,\mathfrak{A}_{\sigma(k+q-1)},b_{k+q})+\\
 &Q(\mathfrak{A}_{\sigma(k)},\cdots,\mathfrak{A}_{\sigma(k+q-1)},a_{k+q})\wedge\alpha^n_{\mathfrak{g}}(b_{k+q}),\alpha^n_{\mathfrak{g}}(\mathfrak{A}_{k+q+1}),
 \cdots,\alpha^n_{\mathfrak{g}}(\mathfrak{A}_{p+q}),\alpha^n_{\mathfrak{g}}(c)),
\end{align*}
for any $\mathfrak{A}_i=a_i\wedge b_i\in \wedge^2 \mathfrak{g},i=1,\cdots,p+q, c\in \mathfrak{g},  $  where $\alpha^n_{\mathfrak{g}}( \mathfrak{A}_i)=\alpha^n_{\mathfrak{g}}(a_i)\wedge\alpha^n_{\mathfrak{g}}(b_i),$
$\mathbb{S}_{(i,n-i)}$ denotes the set
of $(i,n-i)$-unshuffles, i.e.,  the permutation $\sigma\in\mathbb{S}_{n}$  satisfies $ \sigma(1)<\cdots<\sigma(i)$ and  $ \sigma(i+1)<\cdots<\sigma(n)$.

It is important to note that the coboundary operator for the   cohomology of the
Hom-Lie triple system $(\mathfrak{g},[-,-,-]_\mathfrak{g}, \alpha_\mathfrak{g})$ with coefficients in itself is simply given by
$$\delta f=(-1)^{n-1}[\pi,f]_{\mathrm{Hlts}}, \forall f\in C^n_{\mathrm{H}}(\mathfrak{g},\mathfrak{g}),$$
where  $\pi$ corresponds to the Hom-Lie triple system bracket  $[-,-,-]_\mathfrak{g}$.

\begin{prop}
Assume that  $\pi\in C^2_\mathrm{H}(\mathfrak{g},\mathfrak{g})$
satisfies $\circlearrowleft_{a,b,c}\pi(a,b,c)=0, \forall a,b,c\in \mathfrak{g}$.
Then $(\mathfrak{g},\pi,\alpha_\mathfrak{g})$ is a Hom-Lie triple system  if and only if $\pi$ is a Maurer-Cartan element of the graded Lie
algebra $(C^*_\mathrm{H}(\mathfrak{g},\mathfrak{g}),[-,-]_{\mathrm{Hlts}})$,  i.e., $[\pi,\pi]_{\mathrm{Hlts}}=0$.
\end{prop}

\begin{proof}
Let $\pi\in C^2_\mathrm{H}(\mathfrak{g},\mathfrak{g})$, we have
 \begin{align*}
\pi\diamond\pi(\mathfrak{A}_1,\mathfrak{A}_2,c)&=-\pi(\alpha_{\mathfrak{g}}(\mathfrak{A}_1),\pi(\mathfrak{A}_2,c))+\pi(\alpha_{\mathfrak{g}}(a_2)\wedge\pi(\mathfrak{A}_1,b_2),\alpha_{\mathfrak{g}}(c))+\\
&\pi(\pi(\mathfrak{A}_1,a_2)\wedge\alpha_{\mathfrak{g}}(b_2),\alpha_{\mathfrak{g}}(c))+\pi(\alpha_{\mathfrak{g}}(\mathfrak{A}_2),\pi(\mathfrak{A}_1,c)),
\end{align*}
which implies that $\pi$ defines a Hom-Lie triple system structure on $\mathfrak{g}$ if and only if  $[\pi,\pi]_{\mathrm{Hlts}}=0$.
\end{proof}

Now, we recall some known results on bidifferential graded Lie algebras given in \cite{Liu23}.

\begin{defn} \cite{Liu23}
A differential graded Lie algebra is a triple $(\mathfrak{g}=\oplus_{i=0}^{+\infty} \mathfrak{g}^i,[-,-],d)$  consisting of a
graded Lie algebra  $(\mathfrak{g}=\oplus_{i=0}^{+\infty} \mathfrak{g}^i,[-,-])$  and a differential  $d:\mathfrak{g}\rightarrow\mathfrak{g}$  satisfying  $d(\mathfrak{g}_k)\subseteq\mathfrak{g}_{k+1}$  and
$$d([a,b])=[da,b]+(-1)^k[a,db],\forall a,b\in \mathfrak{g}_k.$$
In particular, any graded Lie algebra is a differential graded Lie algebra with  $d=0.$
Furthermore, an element $\pi\in \mathfrak{g}_1$ is called a Maurer-Cartan element of the differential graded Lie algebra
$(\mathfrak{g}=\oplus_{i=0}^{+\infty} \mathfrak{g}^i,[-,-],d)$ if it satisfies $d\pi+\frac{1}{2}[\pi,\pi]=0.$
\end{defn}

\begin{defn} \cite{Liu23}
Let $(\mathfrak{g}=\oplus_{i=0}^{+\infty} \mathfrak{g}^i,[-,-],d_1)$ and $(\mathfrak{g}=\oplus_{i=0}^{+\infty} \mathfrak{g}^i,[-,-],d_2)$ be two differential graded Lie algebras. The quadruple
$(\mathfrak{g}=\oplus_{i=0}^{+\infty} \mathfrak{g}^i,[-,-],d_1,d_2)$ is called a bidifferential graded Lie algebra if $d_1$ and $d_2$ satisfy
$$d_1\circ d_2+d_2\circ d_1=0.$$
\end{defn}

\begin{defn} \cite{Liu23}
 Let $(\mathfrak{g}=\oplus_{i=0}^{+\infty} \mathfrak{g}^i,[-,-],d_1,d_2)$ be a bidifferential graded Lie algebra. A pair $(\pi_1,\pi_2)\in \mathfrak{g}_1\oplus\mathfrak{g}_1$ is
said to be a Maurer-Cartan element of the bidifferential graded Lie algebra $(\mathfrak{g},[-,-],d_1,d_2)$, if $\pi_1$ and $\pi_2$
are Maurer-Cartan elements of the differential graded Lie algebras $(\mathfrak{g},[-,-],d_1)$  and $(\mathfrak{g},[-,-],d_2)$  respectively,
and
$$d_1\pi_2+ d_2\pi_1+[\pi_1,\pi_2]=0.$$
\end{defn}

\begin{prop}\cite{Liu23} \label{prop:2.8}
Let $(\mathfrak{g},[-,-],d_1,d_2)$  be a bidifferential graded Lie algebra and $(\pi_1,\pi_2)$ be its
Maurer-Cartan element. Then $(\mathfrak{g},[-,-],\widetilde{d}_1,\widetilde{d}_2)$ is a bidifferential graded Lie algebra, where $\widetilde{d}_1$ and  $\widetilde{d}_2$  are
given by
 $$\widetilde{d}_1\mu=d_1\mu+[\pi_1,\mu], \widetilde{d}_2\mu=d_1\mu+[\pi_2,\mu], \forall \mu\in \mathfrak{g}.$$
Moreover, for all $\widetilde{\pi}_1,\widetilde{\pi}_2\in \mathfrak{g}_1$, $(\pi_1+\widetilde{\pi}_1,\pi_2+\widetilde{\pi}_2)$ is a Maurer-Cartan element of the bidifferential graded
Lie algebra $(\mathfrak{g},[-,-],d_1,d_2)$ if and only if the pair $(\widetilde{\pi}_1,\widetilde{\pi}_2)$ is a Maurer-Cartan element of the bidifferential
graded Lie algebra $(\mathfrak{g},[-,-],d_1,d_2)$.
\end{prop}

Next we introduce the concept of a compatible Hom-Lie triple system and provide its Maurer-Cartan characterization. We also define the representation of a compatible Hom-Lie triple system.

\begin{defn}
A compatible Hom-Lie triple system is a  quadruple $(\mathfrak{g},[-,-,-]^1_\mathfrak{g},[-,-,-]^2_\mathfrak{g},$ $\alpha_\mathfrak{g})$ in which $(\mathfrak{g},[-,-,-]^1_\mathfrak{g}, \alpha_\mathfrak{g})$  and
$(\mathfrak{g}, [-,-,-]^2_\mathfrak{g},\alpha_\mathfrak{g})$ are Hom-Lie triple system satisfying
\begin{align}
&[\alpha_\mathfrak{g}(a), \alpha_\mathfrak{g}(b), [x, y, z]^2_\mathfrak{g}]^1_\mathfrak{g}+[\alpha_\mathfrak{g}(a), \alpha_\mathfrak{g}(b), [x, y, z]^1_\mathfrak{g}]^2_\mathfrak{g}\nonumber\\
=&[[a, b, x]^1_\mathfrak{g}, \alpha_\mathfrak{g}(y), \alpha_\mathfrak{g}(z)]^2_\mathfrak{g}+ [[a, b, x]^2_\mathfrak{g}, \alpha_\mathfrak{g}(y), \alpha_\mathfrak{g}(z)]^1_\mathfrak{g}+[\alpha_\mathfrak{g}(x),  [a, b, y]^1_\mathfrak{g}, \alpha_\mathfrak{g}(z)]^2_\mathfrak{g}\nonumber\\
&+[\alpha_\mathfrak{g}(x),  [a, b, y]^2_\mathfrak{g}, \alpha_\mathfrak{g}(z)]^1_\mathfrak{g}+ [\alpha_\mathfrak{g}(x),  \alpha_\mathfrak{g}(y), [a, b, z]^1_\mathfrak{g}]^2_\mathfrak{g}+ [\alpha_\mathfrak{g}(x),  \alpha_\mathfrak{g}(y), [a, b, z]^2_\mathfrak{g}]^1_\mathfrak{g}\label{2.7}
\end{align}
\end{defn}

A homomorphism between two compatible Hom-Lie triple systems  $(\mathfrak{g}_1, [-, -, -]^1_{\mathfrak{g}_1},$ $[-, -, -]^2_{\mathfrak{g}_1},\alpha_1)$ and $(\mathfrak{g}_2, [-, -, -]^1_{\mathfrak{g}_2},[-, -, -]^2_{\mathfrak{g}_2},\alpha_2)$ is a linear map $\varphi: \mathfrak{g}_1\rightarrow \mathfrak{g}_2$ satisfying
$$\varphi(\alpha_1(x))=\alpha_2(\varphi(x)), ~~\varphi([x, y, z]^i_{\mathfrak{g}_1})=[\varphi(x), \varphi(y),\varphi(z)]^i_{\mathfrak{g}_2},~~ {\text {for}} ~~i=1,2 ~~{\text {and }} ~x,y,z\in \mathfrak{g}_1.$$
Furthermore, if $\varphi$ is  nondegenerate,  then $\varphi$ is called an isomorphism from $\mathfrak{g}_1$ to $\mathfrak{g}_2$.

\begin{exam}
 Consider a 2-dimensional Hom-Lie triple system  $\mathfrak{g}$ with basis $\varepsilon_1,\varepsilon_2$ defined by
$$[\varepsilon_1,\varepsilon_2,\varepsilon_2]^1_\mathfrak{g}=\varepsilon_1, [\varepsilon_1,\varepsilon_2,\varepsilon_2]^2_\mathfrak{g}=\varepsilon_2, \alpha_\mathfrak{g}(\varepsilon_1)=\varepsilon_1, \alpha_\mathfrak{g}(\varepsilon_2)=-\varepsilon_2.$$
Then $(\mathfrak{g},[-,-,-]^1_\mathfrak{g},[-,-,-]^2_\mathfrak{g}, \alpha_\mathfrak{g})$  is a compatible Hom-Lie triple system.
\end{exam}

\begin{exam}
Let $(\mathfrak{g},[-,-]^1_\mathfrak{g},[-,-]^2_\mathfrak{g}, \alpha_\mathfrak{g})$ be a compatible Hom-Lie algebras (see \cite{Das23}). We define
 $$[-,-,-]^i_\mathfrak{g}:\wedge^2 \mathfrak{g} \otimes \mathfrak{g}\rightarrow \mathfrak{g}~~~ \text{by}~~~ [a,b,c]^i_\mathfrak{g}=[[a,b]^j_\mathfrak{g},\alpha_\mathfrak{g}(c)]^i_\mathfrak{g}, \a,b,c\in \mathfrak{g},i,j=1,2.$$
Then $(\mathfrak{g},[-,-,-]^1_\mathfrak{g},[-,-,-]^2_\mathfrak{g}, \alpha_\mathfrak{g})$  is a compatible Hom-Lie triple system.
\end{exam}

Let $(\mathfrak{g}, [-, -, -]_\mathfrak{g},\alpha_\mathfrak{g})$ be a Hom-Lie triple system.
A Nijenhuis operator on $(\mathfrak{g}, [-, -, -]_\mathfrak{g},\alpha_\mathfrak{g})$  is  a linear map $N:\mathfrak{g}\rightarrow\mathfrak{g}$  satisfying
\begin{align*}
 N\circ \alpha_\mathfrak{g}=&\alpha_\mathfrak{g}\circ N,\\
 [Nx,Ny,Nz]_\mathfrak{g}=&N([x,Ny,Nz]_\mathfrak{g}+[Nx,y,Nz]_\mathfrak{g}+[Nx,Ny,z]_\mathfrak{g})\\
 &-N^2([Nx,y,z]_\mathfrak{g}+[x,Ny,z]_\mathfrak{g}+[x,y,Nz]_\mathfrak{g})+N^3[x,y,z]_\mathfrak{g},
\end{align*}
 for all $x,y,z\in \mathfrak{g}$.
 Nijenhuis operators are useful to study  2-order  linear deformation of a Hom-Lie triple system \cite{Hou20}.

\begin{exam}
Let $N$ be a Nijenhuis operator on a Hom-Lie triple system $(\mathfrak{g}, [-, -, -]_\mathfrak{g},\alpha_\mathfrak{g})$. Then  $(\mathfrak{g}, [-, -, -]_N,\alpha_\mathfrak{g})$ is a  Hom-Lie triple system, where
\begin{align*}
 [x,y,z]_N=&[x,Ny,Nz]_\mathfrak{g}+[Nx,y,Nz]_\mathfrak{g}+[Nx,Ny,z]_\mathfrak{g}\\
 &-N([Nx,y,z]_\mathfrak{g}+[x,Ny,z]_\mathfrak{g}+[x,y,Nz]_\mathfrak{g})+N^2[x,y,z]_\mathfrak{g}.
\end{align*}
Furthermore,  the quadruple $(\mathfrak{g}, [-, -, -]_\mathfrak{g}, [-, -, -]_N,\alpha_\mathfrak{g})$ is a
compatible Hom-Lie triple system.
\end{exam}

The proof of the following proposition is straightforward.
\begin{prop}
 A  quadruple  $(\mathfrak{g},[-,-,-]^1_\mathfrak{g},[-,-,-]^2_\mathfrak{g}, \alpha_\mathfrak{g})$  is a compatible Hom-Lie triple system if and only if $([-,-,-]^1_\mathfrak{g}, \alpha_\mathfrak{g})$ and $([-,-,-]^2_\mathfrak{g}, \alpha_\mathfrak{g})$ are Hom-Lie triple system
structures on L such that for any $k_1,k_2\in \mathbb{K}$, the new
trilinear operation
$$\{-,-,-\}_\mathfrak{g}=k_1[-,-,-]^1_\mathfrak{g}+k_2[-,-,-]^2_\mathfrak{g}$$
define the  Hom-Lie triple system structure on Hom-vector space $(\mathfrak{g},  \alpha_\mathfrak{g})$.
\end{prop}

\begin{prop} \label{prop:2.14}
Assume that  $\pi_i\in C^2_\mathrm{H}(\mathfrak{g},\mathfrak{g})$
satisfies $\circlearrowleft_{a,b,c}\pi_i(a,b,c)=0, $ for $i=1,2,a,b,c\in \mathfrak{g}$.
Then  $(\mathfrak{g} , \pi_1, \pi_2, \alpha_\mathfrak{g})$ is a compatible Hom-Lie triple system if and only if $(\pi_1,\pi_2)$ is a Maurer-Cartan element of the   bidifferential
graded Lie algebra $(C^*_\mathrm{H}(\mathfrak{g},\mathfrak{g}),[-,-,$ $-]_{\mathrm{Hlts}},d_1=0,d_2=0)$.
\end{prop}

\begin{proof}
Let $\pi_1,\pi_2\in C^2_\mathrm{H}(\mathfrak{g},\mathfrak{g})$, then, we have
 \begin{align*}
\pi_1 ~\text{defines a Hom-Lie triple system structure} &\Leftrightarrow [\pi_1,\pi_1]_{\mathrm{Hlts}}=0,\\
\pi_2 ~\text{defines a Hom-Lie triple system structure} &\Leftrightarrow [\pi_2,\pi_2]_{\mathrm{Hlts}}=0,\\
\text{compatibility equation ~\eqref{2.7}} &\Leftrightarrow [\pi_1,\pi_2]_{\mathrm{Hlts}}=0.
\end{align*}
Therefore, $(\mathfrak{g} , \pi_1, \pi_2, \alpha_\mathfrak{g})$ is a compatible Hom-Lie triple system if and only if $(\pi_1,\pi_2)$ is a Maurer-Cartan element of the   bidifferential
graded Lie algebra $(C^*_\mathrm{H}(\mathfrak{g},\mathfrak{g}),[-,-,$ $-]_{\mathrm{Hlts}},d_1=0,d_2=0)$.
\end{proof}

So, from Proposition \ref{prop:2.8} and Proposition \ref{prop:2.14}, we get the following conclusion.

\begin{prop}
Let $(\mathfrak{g} , \pi_1, \pi_2, \alpha_\mathfrak{g})$ be a compatible Hom-Lie triple system.

Then, \text{(i)}   $(C^*_{\mathrm{H}}(\mathfrak{g},\mathfrak{g}),[-,-]_{\mathrm{Hlts}},d_{\pi_1},d_{\pi_2})$  is  a bidifferential graded Lie algebra, where
 $d_{\pi_1} Q=[\pi_1,Q]_{\mathrm{Hlts}},~~d_{\pi_2} Q=[\pi_2,Q]_{\mathrm{Hlts}},~~\forall Q\in C^p_{\mathrm{H}}(\mathfrak{g},\mathfrak{g}).$

\text{(ii)} For all $\widetilde{\pi}_1,\widetilde{\pi}_2\in C^2_{\mathrm{H}}(\mathfrak{g},\mathfrak{g})$ satisfying $\circlearrowleft_{a,b,c}\widetilde{\pi}_i(a,b,c)=0, \forall a,b,c\in \mathfrak{g},i=1,2,$ the quadruple $(\mathfrak{g} , \pi_1+\widetilde{\pi}_1,  \pi_2+\widetilde{\pi}_2, \alpha_\mathfrak{g})$ is a compatible Hom-Lie triple system if and only if  $(\widetilde{\pi}_1,\widetilde{\pi}_2)$ is a Maurer-Cartan element of the bidifferential graded
Lie algebra $(C^*_{\mathrm{H}}(\mathfrak{g},\mathfrak{g}),[-,$ $-]_{\mathrm{Hlts}},d_{\pi_1},d_{\pi_2})$ .
\end{prop}

Next, we define representations of a compatible Hom-Lie triple system.

\begin{defn}
A representation of the compatible Hom-Lie triple system  $(\mathfrak{g},[-,-,-]^1_\mathfrak{g},$ $[-,-,-]^2_\mathfrak{g},$ $\alpha_\mathfrak{g})$  consists of a
quadruple $(V,\theta_1,\theta_2, \beta)$  such that

(i) $(V,\theta_1,  \beta)$ is a representation of the Hom-Lie triple system  $(\mathfrak{g},[-,-,-]^1_\mathfrak{g}, $ $\alpha_\mathfrak{g})$;

(ii) $(V,\theta_2,  \beta)$ is a representation of the Hom-Lie triple system  $(\mathfrak{g},[-,-,-]^2_\mathfrak{g}, $ $\alpha_\mathfrak{g})$;

(iii)  for any $a_1\in \mathfrak{g}$, the following compatibility equations hold:
\begin{align}
&\theta_1(\alpha_\mathfrak{g}(a_1),   [a_2, a_3, a_4]^2_\mathfrak{g})\circ\beta+\theta_2(\alpha_\mathfrak{g}(a_1),   [a_2, a_3, a_4]^1_\mathfrak{g})\circ\beta+\theta_1(\alpha_\mathfrak{g}(a_2), \alpha_\mathfrak{g}(a_4))\theta_2(a_1,  a_3)\nonumber\\
&+\theta_2(\alpha_\mathfrak{g}(a_2),  \alpha_\mathfrak{g}(a_4))\theta_1(a_1,  a_3)-\theta_1(\alpha_\mathfrak{g}(a_3),    \alpha_\mathfrak{g}(a_4))\theta_2(a_1,  a_2) -\theta_2(\alpha_\mathfrak{g}(a_3), \alpha_\mathfrak{g}(a_4))\theta_1(a_1,  a_2)\nonumber\\
&-D_1(\alpha_\mathfrak{g}(a_2),  \alpha_\mathfrak{g}(a_3))\theta_2(a_1,  a_4)-D_2(\alpha_\mathfrak{g}(a_2),    \alpha_\mathfrak{g}(a_3))\theta_1(a_1,  a_4)=0,\label{2.8}\\
&D_1(\alpha_\mathfrak{g}(a_3),   [a_1, a_2, a_4]^2_\mathfrak{g})\circ\beta+D_2(\alpha_\mathfrak{g}(a_3),   [a_1, a_2, a_4]^1_\mathfrak{g})\circ\beta+
D_1(  [a_1, a_2, a_3]^2_\mathfrak{g},\alpha_\mathfrak{g}(a_4))\circ\beta+\nonumber\\
&D_2([a_1, a_2, a_3]^1_\mathfrak{g},\alpha_\mathfrak{g}(a_4))\circ\beta-D_1(\alpha_\mathfrak{g}(a_1), \alpha_\mathfrak{g}(a_2))D_2(a_3,  a_4)+D_2(\alpha_\mathfrak{g}(a_3),  \alpha_\mathfrak{g}(a_4))D_1(a_1,  a_2)+\nonumber\\
&D_1(\alpha_\mathfrak{g}(a_3),    \alpha_\mathfrak{g}(a_4))D_2(a_1,  a_2)-D_2(\alpha_\mathfrak{g}(a_1), \alpha_\mathfrak{g}(a_2))D_1(a_3,  a_4)=0,\label{2.9}
\end{align}
where $D_k(a_i,a_j)=\theta_k(a_j,a_i)-\theta_k(a_i,a_j), k=1,2$.
\end{defn}

It follows that any compatible Hom-Lie triple system $(\mathfrak{g},[-,-,-]^1_\mathfrak{g},$ $[-,-,-]^2_\mathfrak{g}, \alpha_\mathfrak{g})$ is a representation of itself, where
$V=\mathfrak{g},\theta_1(a,b)c=[c,a,b]^1_\mathfrak{g},\theta_2(a,b)c=[c,a,b]^2_\mathfrak{g}$ and $\beta=\alpha_\mathfrak{g}$. This is called the adjoint representation.

\begin{remark} \label{remark:2.8}
Let $(V,\theta_1,\theta_2, \beta)$   be a representation of a compatible Hom-Lie triple system  $(\mathfrak{g},[-,-,-]^1_\mathfrak{g},$ $[-,-,-]^2_\mathfrak{g},$ $\alpha_\mathfrak{g})$.
Then $(V, k_1\theta_1+ k_2\theta_2, \beta)$ is a representation of the Hom-Lie triple system $(\mathfrak{g}, k_1[-,-,-]^1_\mathfrak{g}+k_2[-,-,-]^2_\mathfrak{g}, \alpha_\mathfrak{g})$.
\end{remark}

It is similar to the standard case that the following proposition is proved.

\begin{prop} \label{prop:semi-product}
Let $(V,\theta_1,\theta_2, \beta)$   be a representation of a compatible Hom-Lie triple system  $(\mathfrak{g},[-,-,-]^1_\mathfrak{g},$ $[-,-,-]^2_\mathfrak{g},$ $\alpha_\mathfrak{g})$.
Define two trilinear operations $[-,-,-]^1_{\ltimes},[-,-,-]^2_{\ltimes}:\wedge^2(\mathfrak{g}\oplus V)\otimes(\mathfrak{g}\oplus V)\rightarrow  \mathfrak{g}\oplus V$
and a linear map $\alpha_\mathfrak{g}\oplus\beta: \mathfrak{g}\oplus V \rightarrow  \mathfrak{g}\oplus V$ by
\begin{align*}
[(a,u),(b,v),(c,w)]^1_{\ltimes}&=([a,b,c]^1_\mathfrak{g},D_1(a,b)w+\theta_1(b,c)u-\theta_1(a,c)v),\\
[(a,u),(b,v),(c,w)]^2_{\ltimes}&=([a,b,c]^2_\mathfrak{g},D_2(a,b)w+\theta_2(b,c)u-\theta_2(a,c)v),\\
\alpha_\mathfrak{g}\oplus\beta(a,u)&=(\alpha_\mathfrak{g}(a),\beta(u)),
\end{align*}
for any $(a,u),(b,v),(c,w)\in \mathfrak{g}\oplus V$.  Then $(\mathfrak{g}\oplus V,[-,-,-]^1_{\ltimes},[-,-,-]^2_{\ltimes}, \alpha_\mathfrak{g}\oplus\beta)$ is a compatible Hom-Lie triple system.
This compatible Hom-Lie triple system is called the semi-product
compatible Hom-Lie triple system and denoted by  $\mathfrak{g}\ltimes^1_2 V$.
\end{prop}

\section{   Cohomology of compatible   Hom-Lie triple systems }\label{sec: Cohomologies}
\def\theequation{\arabic{section}.\arabic{equation}}
\setcounter{equation} {0}

In this section, we introduce the cohomology of a compatible Hom-Lie triple system with coefficients in a representation.

Next,   we always assume that $(V,\theta_1,\theta_2, \beta)$   is a representation of a compatible Hom-Lie triple system  $(\mathfrak{g},[-,-,-]^1_\mathfrak{g},$ $[-,-,-]^2_\mathfrak{g},$ $\alpha_\mathfrak{g})$.

Let $\delta_1: \mathcal{C}_{\mathrm{HLts}}^{n}(\mathfrak{g}, V)\rightarrow \mathcal{C}_{\mathrm{HLts}}^{n+1}( \mathfrak{g}, V)$(resp. $\delta_2: \mathcal{C}_{\mathrm{HLts}}^{n}(\mathfrak{g}, V)\rightarrow \mathcal{C}_{\mathrm{HLts}}^{n+1}( \mathfrak{g}, V)$), for $n\geq 1$, be
the coboundary operator for the cohomology of the Hom-Lie triple system  $(\mathfrak{g},[-,-,-]^1_\mathfrak{g},\alpha_\mathfrak{g})$ with coefficients in the representation $(V,\theta_1,  \beta)$  (resp. of the Hom-Lie triple system $(\mathfrak{g},[-,-,-]^2_\mathfrak{g},\alpha_\mathfrak{g})$  with coefficients in the representation  $(V,\theta_2, \beta)$).  Then, we have
$(\delta_1)^2=0$  and $(\delta_2)^2=0$.

Next we consider the semidirect product compatible Hom-Lie triple system structure on $\mathfrak{g}\oplus  V$  given in Proposition  \ref{prop:semi-product}.
Notice that any   map  $f\in \mathcal{C}^{n-1}_{\mathrm{HLts}}(\mathfrak{g},V)$  can be lifted to a   map $\widetilde{f}\in C^{n-1}_{\mathrm{H}}(\mathfrak{g}\oplus V,\mathfrak{g}\oplus V)$  by
$$\widetilde{f}((\mathfrak{A}_1,\mathfrak{V}_1),\cdots,(\mathfrak{A}_{n-1},\mathfrak{V}_{n-1}),(c,u))=(0,f(\mathfrak{A}_1,\cdots,\mathfrak{A}_{n-1},c)),$$
for any $\mathfrak{A}_i=a_i\wedge b_i\in \wedge^2\mathfrak{g},\mathfrak{V}_i=u_i\wedge v_i\in \wedge^2 V, c\in \mathfrak{g},u\in V$. Then $f=0$ if and only if $\widetilde{f}=0$.

 We denote by  $\pi_1,\pi_2\in C^{2}_{\mathrm{H}}(\mathfrak{g}\oplus V,\mathfrak{g}\oplus V)$   the elements
corresponding to the Hom-Lie triple system structures $[-,-,-]^1_{\ltimes}$
and $[-,-,-]^2_{\ltimes}$
on $\mathfrak{g}\oplus V$, respectively. Let
$$\widetilde{\delta}_1: C_{\mathrm{H}}^{n}(\mathfrak{g}\oplus V, \mathfrak{g}\oplus V)\rightarrow C_{\mathrm{H}}^{n+1}(\mathfrak{g}\oplus V, \mathfrak{g}\oplus V),$$
$$\widetilde{\delta}_2: C_{\mathrm{H}}^{n}(\mathfrak{g}\oplus V, \mathfrak{g}\oplus V)\rightarrow C_{\mathrm{H}}^{n+1}(\mathfrak{g}\oplus V, \mathfrak{g}\oplus V)$$
denote respectively the coboundary operator for the cohomology of the Hom-Lie triple system
$(\mathfrak{g}\oplus V,\pi_1, \alpha_\mathfrak{g}\oplus\beta)$ (resp. of the Hom-Lie triple system  $(\mathfrak{g}\oplus V,\pi_2, \alpha_\mathfrak{g}\oplus\beta)$ ) with coefficients in itself.
Further,  for any  $f\in \mathcal{C}^{n-1}_{\mathrm{HLts}}(\mathfrak{g},V)$, we have
$\widetilde{\delta_1f}=\widetilde{\delta}_1(\widetilde{f})=(-1)^{n-1}[\pi_1,\widetilde{f}]_\mathrm{Hlts},\widetilde{\delta_2f}=\widetilde{\delta}_2(\widetilde{f})=(-1)^{n-1}[\pi_2,\widetilde{f}]_\mathrm{Hlts}.$

Then, we have the following conclusion.

\begin{lemma} \label{lemma:3.1}
  The coboundary operators $\delta_1$  and  $\delta_2$  satisfy the following compatibility
  $$ \delta_1\circ\delta_2+\delta_2\circ\delta_1=0.$$
\end{lemma}

\begin{proof}
For any  $f\in \mathcal{C}^{n-1}_{\mathrm{HLts}}(\mathfrak{g},V)$, we have
\begin{align*}
&  \widetilde{(\delta_1\circ\delta_2+\delta_2\circ\delta_1)(f)}\\
&= \widetilde{\delta_1(\delta_2f)}+ \widetilde{\delta_2(\delta_1f)}\\
&= (-1)^n[\pi_1,\widetilde{\delta_2f}]_\mathrm{H}+ (-1)^n[\pi_2,\widetilde{\delta_1f}]_\mathrm{H}\\
&= -[\pi_1,[\pi_2,\widetilde{f}]_\mathrm{H}]_\mathrm{H}-  [\pi_2,[\pi_1,\widetilde{f}]_\mathrm{H}]_\mathrm{H}\\
&= -[[\pi_1,\pi_2]_\mathrm{H},\widetilde{f}]_\mathrm{H}+[\pi_2,[\pi_1,\widetilde{f}]_\mathrm{H}]_\mathrm{H}-  [\pi_2,[\pi_1,\widetilde{f}]_\mathrm{H}]_\mathrm{H}\\
&=0.
\end{align*}
 Hence the result follows.
\end{proof}

Next, we define  the cohomology of a compatible Hom-Lie triple system $(\mathfrak{g},[-,-,-]^1_\mathfrak{g},$ $[-,-,-]^2_\mathfrak{g},$ $\alpha_\mathfrak{g})$ with coefficients in
a represenation $(V,\theta_1,\theta_2, \beta)$.
 For each $n\geq 1$, we define the set of  $n$-cochains  by
\begin{align*}
\mathcal{C}_{\mathrm{cHLts}}^{n}( \mathfrak{g},  V):&=\overbrace{\mathcal{C}_{\mathrm{HLts}}^{n}( \mathfrak{g},  V)\oplus\cdots\oplus\mathcal{C}_{\mathrm{HLts}}^{n}( \mathfrak{g},  V)}^n
\end{align*}
 Define a  linear  map  $\delta_c:\mathcal{C}_{\mathrm{cHLts}}^{n}( \mathfrak{g},  V)\rightarrow \mathcal{C}_{\mathrm{cHLts}}^{n+1}( \mathfrak{g},  V)$,  for $n\geq 1$ by
\begin{align*}
&\text{if}~~f\in \mathcal{C}_{\mathrm{cHLts}}^{1}( \mathfrak{g},  V),~~~~\delta_cf=(\delta_1f,\delta_2f), ~~\\
&\text{if}~~(f_1,\cdots,f_n)\in \mathcal{C}_{\mathrm{cHLts}}^{n}( \mathfrak{g},  V),~~ n\geq2, ~~2\leq i\leq n-1,\\
&~~~~~~~~~~ \delta_c(f_1,\cdots,f_n)=(\delta_1f_1,\delta_1f_2+\delta_2f_{1},\cdots,\underbrace{\delta_1f_i+\delta_2f_{i-1}}_{\text{i-th position}},\cdots,\delta_2f_n).
\end{align*}
Then, we have the following proposition.

\begin{prop}\label{prop:1-cocycle}
 The pair $(\mathcal{C}_{\mathrm{cHLts}}^{*}( \mathfrak{g},  V), \delta_c)$ is a cochain complex. So   $ (\delta_c)^2=0.$
\end{prop}

\begin{proof}
For any $f\in \mathcal{C}_{\mathrm{cHLts}}^{1}( \mathfrak{g},  V)$, by Lemma \ref{lemma:3.1}, we have
$$(\delta_c)^2f=\delta_c(\delta_1f,\delta_2f)=(\delta_1\delta_1f_1,\delta_1\delta_2f+\delta_2\delta_1f,\delta_2\delta_2f)=0.$$
Given any $(f_1,\cdots,f_n)\in \mathcal{C}_{\mathrm{cHLts}}^{n}( \mathfrak{g},  V)$, with  $n\geq 2$, we have
\begin{align*}
(\delta_c)^2(f_1,\cdots,f_n)=&\delta_c(\delta_1f_1,\cdots,\delta_1f_i+\delta_2f_{i-1},\cdots,\delta_2f_n)\\
=&(\delta_1\delta_1f_1,\delta_1\delta_1f_2+\delta_1\delta_2f_1+\delta_2\delta_1f_1,\cdots,\\
&\underbrace{\delta_2\delta_2f_{i-2}+\delta_2\delta_1f_{i-1}+\delta_1\delta_2f_{i-1}+\delta_1\delta_1f_{i}}_{3\leq i\leq n-1},\cdots,\\
&\delta_2\delta_2f_{n-1}+\delta_2\delta_1f_{n}+\delta_1\delta_2f_{n}, \delta_2\delta_2f_{n})\\
=&0.
\end{align*}
Thus, $(\mathcal{C}_{\mathrm{cHLts}}^{*}( \mathfrak{g},  V), \delta_c)$ is a cochain complex.
\end{proof}

Denote the set of $n$-cocycles by $\mathcal{Z}^n_{\mathrm{cHLts}}(\mathfrak{g},V):=\{f\in \mathcal{C}_{\mathrm{cHLts}}^{n}( \mathfrak{g},  V)~|~\delta_c f=0\}$, the set of $n$-coboundaries by $\mathcal{B}^n_{\mathrm{cHLts}}(\mathfrak{g},V):=\{\delta_c f ~|~f\in \mathcal{C}_{\mathrm{cHLts}}^{n-1}( \mathfrak{g},  V)\}$, and $n$-th cohomology group by
$$\mathcal{H}^n_{\mathrm{cHLts}}(\mathfrak{g},V)=\frac{\mathcal{Z}^n_{\mathrm{cHLts}}(\mathfrak{g},V)}{\mathcal{B}^n_{\mathrm{cHLts}}(\mathfrak{g},V)}, n\geq 1.$$

Let $(\mathfrak{g},[-,-,-]^1_\mathfrak{g},$ $[-,-,-]^2_\mathfrak{g},$ $\alpha_\mathfrak{g})$  be a compatible Hom-Lie triple system and $(V,\theta_1,\theta_2, \beta)$ be a representation of it.
Then by Remark  \ref{remark:2.8},  $V_+=(V, \theta_1+ \theta_2, \beta)$ is a representation of the Hom-Lie triple system $\mathfrak{g}_+=(\mathfrak{g}, [-,-,-]^1_\mathfrak{g}+[-,-,-]^2_\mathfrak{g}, \alpha_\mathfrak{g})$.
Consider the cochain complex $(\mathcal{C}_{\mathrm{cHLts}}^{*}( \mathfrak{g},  V), \delta_c)$ of the compatible Hom-Lie triple system  $\mathfrak{g}$ with coefficients in the representation $V$, and the cochain complex  $(\mathcal{C}_{\mathrm{HLts}}^{*}( \mathfrak{g}_+,  V_+), \delta)$ of the
Hom-Lie triple system $\mathfrak{g}_+$ with coefficients in $V_+$.

 For each $n\geq 1$, define a  map $\Phi_n: \mathcal{C}_{\mathrm{cHLts}}^{n}( \mathfrak{g},  V)\rightarrow \mathcal{C}_{\mathrm{HLts}}^{*}( \mathfrak{g}_n,  V_+)$ by
$$\Phi_n(f_1,\cdots,f_n)=f_1+\cdots+f_n, \forall ~(f_1,\cdots,f_n)\in \mathcal{C}_{\mathrm{cHLts}}^{n}( \mathfrak{g},  V).$$
Then, we have
\begin{align*}
(\delta\circ\Phi_n)(f_1,\cdots,f_n)&=\delta(f_1+\cdots+f_n)\\
&=\delta_1(f_1,\cdots,f_n)+\delta_2(f_1,\cdots,f_n)\\
&=\Phi_{n+1}(\delta_1f_1, \cdots, \delta_1f_i+\delta_2f_{i-1},\cdots,\delta_2f_n)\\
&=(\Phi_{n+1}\circ\delta_c)(f_1,\cdots,f_n).
\end{align*}
So we have the following theorem.
\begin{theorem}
 The collection $\{\Phi_n\}_{n\geq 1}$
defines a morphism of cochain complexes from $(\mathcal{C}_{\mathrm{cHLts}}^{*}( \mathfrak{g},  V), \delta_c)$
to $(\mathcal{C}_{\mathrm{HLts}}^{*}( \mathfrak{g}_+,  V_+), \delta)$. Hence, it induces a morphism $\mathcal{H}_{\mathrm{cHLts}}^{*}( \mathfrak{g},  V)\rightarrow\mathcal{H}_{\mathrm{HLts}}^{*}( \mathfrak{g}_+,  V_+)$ between corresponding  cohomologies.
\end{theorem}

\section{   Linear deformatons of     compatible Hom-Lie triple systems } \label{sec: deformations}
\def\theequation{\arabic{section}.\arabic{equation}}
\setcounter{equation} {0}

In this section, we study linear deformations  of a compatible Hom-Lie triple system. We introduce Nijenhuis operators that
generate trivial linear deformations.

\begin{defn}
Let $(\mathfrak{g},[-,-,-]^1_\mathfrak{g},$ $[-,-,-]^2_\mathfrak{g},$ $\alpha_\mathfrak{g})$  be a compatible Hom-Lie triple system with  $\pi_1=[-,-,-]^1_\mathfrak{g},\pi_2=[-,-,-]^2_\mathfrak{g}$
   and  $(\mu_1,\mu_2),(\omega_1,\omega_2)\in \mathcal{C}_{\mathrm{cHLts}}^{2}( \mathfrak{g},  \mathfrak{g})$ be two 2-cochains. Define
two trilinear operations on $\mathfrak{g}$ depending on the parameter t as follows:
 \begin{align*}
&[a,b,c]_t^1=\pi_1(a,b,c)+t\mu_1(a,b,c)+t^2\omega_1(a,b,c),\\
& [a,b,c]_t^2=\pi_2(a,b,c)+t\mu_2(a,b,c)+t^2\omega_2(a,b,c), ~\forall ~a,b,c\in\mathfrak{g}.
\end{align*}
If for all $t$, the quadruple  $(\mathfrak{g},[-,-,-]^1_t, [-,-,-]^2_t,\alpha_\mathfrak{g})$ is a compatible Hom-Lie triple system, then we say that $((\mu_1,\mu_2),(\omega_1,\omega_2))$  generates a linear
deformation of the  compatible Hom-Lie triple system $(\mathfrak{g},[-,-,-]^1_\mathfrak{g}, [-,-,-]^2_\mathfrak{g}, \alpha_\mathfrak{g})$.
\end{defn}

The quadruple $(\mathfrak{g},[-,-,-]^1_t, [-,-,-]^2_t,\alpha_\mathfrak{g})$ is a compatible Hom-Lie triple system  is equivalent to saying that
   \begin{align*}
&[\pi_1+t\mu_1+t^2\omega_1,\pi_1+t\mu_1+t^2\omega_1]_{\mathrm{Hlts}}=0,\\
&[\pi_1+t\mu_1+t^2\omega_1,\pi_2+t\mu_2+t^2\omega_2]_{\mathrm{Hlts}}=0, \\
&[\pi_2+t\mu_2+t^2\omega_2,\pi_2+t\mu_2+t^2\omega_2]_{\mathrm{Hlts}}=0,
\end{align*}
which force that
   \begin{align*}
&[\pi_1, \mu_1]_{\mathrm{Hlts}}=0,[\pi_1,\omega_1]_{\mathrm{Hlts}}+[\mu_1,\mu_1]_{\mathrm{Hlts}}=0,[\mu_1,\omega_1]_{\mathrm{Hlts}}=0,[\omega_1,\omega_1]_{\mathrm{Hlts}}=0,[\omega_1,\omega_2]_{\mathrm{Hlts}}=0,\\
&[\pi_1, \mu_2]_{\mathrm{Hlts}}+[\mu_1, \pi_2]_{\mathrm{Hlts}}=0,[\pi_1,\omega_2]_{\mathrm{Hlts}}+[\omega_1,\pi_2]_{\mathrm{Hlts}}+[\mu_1,\mu_2]_{\mathrm{Hlts}}=0,[\mu_1,\omega_2]_{\mathrm{Hlts}}+[\omega_1,\mu_2]_{\mathrm{Hlts}}=0,\\
&[\pi_2, \mu_2]_{\mathrm{Hlts}}=0,[\pi_2,\omega_2]_{\mathrm{Hlts}}+[\mu_2,\mu_2]_{\mathrm{Hlts}}=0,[\mu_2,\omega_2]_{\mathrm{Hlts}}=0,[\omega_2,\omega_2]_{\mathrm{Hlts}}=0.
\end{align*}
Then
   \begin{align*}
\delta_c(\mu_1,\mu_2)=&(\delta_1\mu_1,\delta_1\mu_2+\delta_2\mu_{1},\delta_2\mu_2)\\
=&([\pi_1, \mu_1]_{\mathrm{Hlts}},[\pi_1, \mu_2]_{\mathrm{Hlts}}+[\pi_2, \mu_1]_{\mathrm{Hlts}},[\pi_2, \mu_2]_{\mathrm{Hlts}})\\
=&0.
\end{align*}
Therefore, $(\mu_1,\mu_2)\in \mathcal{C}_{\mathrm{cHLts}}^{2}( \mathfrak{g},  \mathfrak{g})$ is a 2-cocycle in the cohomology of the compatible Hom-Lie triple system  $\mathfrak{g}$ with coefficients in the adjoint representation. Moreover, $(\mu_1,\mu_2)$ is called the infinitesimal of $([-,-,-]^1_t, [-,-,-]^2_t)$.

\begin{defn}
Let $((\mu_1,\mu_2),(\omega_1,\omega_2))$ and $((\mu'_1,\mu'_2),(\omega'_1,\omega'_2))$ generate linear deformations $(\mathfrak{g},[-,-,-]^1_t, [-,-,-]^2_t,\alpha_\mathfrak{g})$ and $(\mathfrak{g},[-,-,-]'^1_t, [-,-,-]'^2_t,\alpha_\mathfrak{g})$
of a compatible Hom-Lie triple system $\mathfrak{g}$.
They are said to be equivalent if there exists a linear map $N:\mathfrak{g}\rightarrow \mathfrak{g}$
such that $N_t=\mathrm{Id}_{\mathfrak{g}}+tN$ satisfying
   \begin{align}
  & \alpha_\mathfrak{g}\circ N_t=N_t\circ\alpha_\mathfrak{g},\label{4.1}\\
  & N_t[a,b,c]^1_t=[N_ta,N_tb,N_tc]'^1_t,~ N_t[a,b,c]^2_t=[N_ta,N_tb,N_tc]'^2_t,~~ \forall a,b,c\in \mathfrak{g}.\label{4.2}
   \end{align}
\end{defn}

We know that Eqs. \eqref{4.1} and \eqref{4.2} are equivalent to
   \begin{align}
  \alpha_\mathfrak{g}\circ N=&N\circ\alpha_\mathfrak{g},\label{4.3}\\
  N\pi_i(a,b,c)+\mu_i(a,b,c)=&\pi_i(Na,b,c)+\pi_i(a,Nb,c)+\pi_i(a,b,Nc)+\mu'_i(a,b,c),\label{4.4}\\
\omega_i(a,b,c)+N\mu_i(a,b,c)=&\pi_i(Na,Nb,c)+\pi_i(Na,b,Nc)+\pi_i(a,Nb,Nc)+\mu'_i(Na,b,c)\nonumber\\
  &+\mu'_i(a,Nb,c)+\mu'_i(a,b,Nc)+\omega'_i(a,b,c),\label{4.5}\\
  N\omega_i(a,b,c)=&\pi_i(Na,Nb,Nc)+\mu'_i(Na,Nb,c)+\mu'_i(Na,b,Nc)+\mu'_i(a,Nb,Nc)\nonumber\\
  &+\omega'_i(Na,b,c)+\omega'_i(a,Nb,c)+\omega'_i(a,b,Nc),\label{4.6}
   \end{align}
      \begin{align}
\mu'_i(Na,Nb,Nc)+\omega'_i(Na,Nb,c)+\omega'_i(Na,b,Nc)+\omega'_i(a,Nb,Nc)&=0,\label{4.7}\\
\omega'_i(Na,Nb,Nc)&=0,\label{4.8}
   \end{align}
   for any $a,b,c\in \mathfrak{g}$ and $i=1,2.$

By  Eq. \eqref{4.4}, we have  $\mu_i-\mu'_i=\delta_i N$, where $N$ is considered as an element in $\mathcal{C}_{\mathrm{cHLts}}^{1}( \mathfrak{g},   \mathfrak{g})$. So $(\mu_1-\mu'_1,\mu_2-\mu'_2)=(\delta_1 N, \delta_2 N)=\delta_c N\in \mathcal{B}_{\mathrm{cHLts}}^{2}( \mathfrak{g},   \mathfrak{g})$.
Hence, their cohomology classes are the same in $ \mathcal{H}_{\mathrm{cHLts}}^{2}( \mathfrak{g},   \mathfrak{g})$.
  To sum up, we
have the following result.
\begin{prop}
 Let   $(\mathfrak{g},[-,-,-]^1_\mathfrak{g},$ $[-,-,-]^2_\mathfrak{g},$ $\alpha_\mathfrak{g})$  be a compatible Hom-Lie triple system.   Then there is a bijection between the set of all
equivalence classes of  infinitesimal of linear deformation of   $\mathfrak{g}$  and the  second  cohomology group $ \mathcal{H}_{\mathrm{cHLts}}^{2}( \mathfrak{g},   \mathfrak{g})$.
\end{prop}

Now   we consider  trivial linear deformations   in order to introduce the notion of a Nijenhuis operator on a     compatible  Hom-Lie triple system.

\begin{defn}
A linear deformation $(\pi_1+t\mu_1+t^2\omega_1, \pi_2+t\mu_2+t^2\omega_2)$  of compatible Hom-Lie triple system $(\mathfrak{g},[-,-,-]^1_\mathfrak{g},$ $[-,-,-]^2_\mathfrak{g},$ $\alpha_\mathfrak{g})$
is said to be trivial if the deformation is equivalent to the undeformed one $(\pi_1, \pi_2)$.
\end{defn}

Therefore, a linear deformation $(\pi_1+t\mu_1+t^2\omega_1, \pi_2+t\mu_2+t^2\omega_2)$  is trivial if and only if there exists a linear map $N:\mathfrak{g}\rightarrow \mathfrak{g}$
 satisfying
   \begin{align}
  \alpha_\mathfrak{g}\circ N=&N\circ\alpha_\mathfrak{g},\label{4.9}\\
  N\pi_i(a,b,c)+\mu_i(a,b,c)=&\pi_i(Na,b,c)+\pi_i(a,Nb,c)+\pi_i(a,b,Nc),\label{4.10}\\
\omega_i(a,b,c)+N\mu_i(a,b,c)=&\pi_i(Na,Nb,c)+\pi_i(Na,b,Nc)+\pi_i(a,Nb,Nc),\label{4.11}\\
  N\omega_i(a,b,c)=&\pi_i(Na,Nb,Nc),\label{4.12}
   \end{align}
  for any $a,b,c\in \mathfrak{g}$ and $i=1,2.$

\begin{defn}
A Nijenhuis operator on a compatible Hom-Lie triple system  $(\mathfrak{g},[-,-,-]^1_\mathfrak{g},$ $[-,-,-]^2_\mathfrak{g},$ $\alpha_\mathfrak{g})$ is a linear map
$N:\mathfrak{g}\rightarrow \mathfrak{g}$ which is a Nijenhuis operator for both the Hom-Lie triple systems $(\mathfrak{g},[-,-,-]^1_\mathfrak{g},$   $\alpha_\mathfrak{g})$  and $(\mathfrak{g}, $ $[-,-,-]^2_\mathfrak{g},$ $\alpha_\mathfrak{g})$ , i.e., $N$
satisfies
   \begin{align*}
  \alpha_\mathfrak{g}\circ N=&N\circ\alpha_\mathfrak{g}, \\
  [Na,Nb,Nc]^i_\mathfrak{g}=&N([Na,Nb,c]^i_\mathfrak{g}+[Na,b,Nc]^i_\mathfrak{g}+[Na,Nb,Nc]^i_\mathfrak{g})\\
  &-N^2([Na,b,c]^i_\mathfrak{g}+[a,Nb,c]^i_\mathfrak{g}+[a,b,Nc]^i_\mathfrak{g})+N^3[a,b,c]^i_\mathfrak{g},~~~\text{for}~ i=1,2.
   \end{align*}
\end{defn}

Explicitly, any trivial linear deformation of  a compatible Hom-Lie algebra can characterize the  Nijenhuis operator.
Next, the converse is given by the next result whose proof is straightforward.

\begin{prop}
 Let $N:\mathfrak{g}\rightarrow \mathfrak{g}$  be a Nijenhuis operator on a compatible Hom-Lie triple system $(\mathfrak{g},[-,-,-]^1_\mathfrak{g},$ $[-,-,-]^2_\mathfrak{g},$ $\alpha_\mathfrak{g})$.
Then $((\mu_1,\mu_2),(\omega_1,\omega_2))$  generates a trivial linear deformation of the compatible Hom-Lie triple system $(\mathfrak{g},[-,-,-]^1_\mathfrak{g},$ $[-,-,-]^2_\mathfrak{g},$ $\alpha_\mathfrak{g})$,
where
   \begin{align*}
 \mu_i(a,b,c)=&[Na,b,c]^i_\mathfrak{g}+[a,Nb,c]^i_\mathfrak{g}+[a,b,Nc]^i_\mathfrak{g}-N[a,b,c]^i_\mathfrak{g}, \\
 \omega_i(a,b,c)=&[Na,Nb,c]^i_\mathfrak{g}+[Na,b,Nc]^i_\mathfrak{g}+[a,Nb,Nc]^i_\mathfrak{g}-N\mu_i(a,b,c), ~\forall a,b,c\in\mathfrak{g}~\text{and}~ i=1,2.
   \end{align*}
\end{prop}

\section{ Abelian extensions of compatible  Hom-Lie triple systems }\label{sec: abelian extensions}
\def\theequation{\arabic{section}.\arabic{equation}}
\setcounter{equation} {0}

In this section,  we study abelian extensions of compatible
Hom-Lie triple systems and give a classification of equivalence classes of abelian extensions.

Let $(\mathfrak{g},[-,-,-]^1_\mathfrak{g}, [-,-,-]^2_\mathfrak{g}, \alpha_\mathfrak{g})$  be a compatible Hom-Lie triple system and $(V, \alpha_V)$ be a Hom-vector space.  Note that $(V, \alpha_V)$ can be considered as a compatible  Hom-Lie triple system with both the Hom-Lie triple system
brackets on $V$ are trivial.

\begin{defn}
Let $(\mathfrak{g},[-,-,-]^1_\mathfrak{g}, [-,-,-]^2_\mathfrak{g}, \alpha_\mathfrak{g})$ and $(V,[-,-,-]^1_V, [-,-,-]^2_V, \alpha_V)$ be two compatible  Hom-Lie triple systems.
An abelian extension of $(\mathfrak{g},[-,-,-]^1_\mathfrak{g}, [-,-,-]^2_\mathfrak{g}, \alpha_\mathfrak{g})$ by  $(V,[-,-,-]^1_V, [-,-,-]^2_V, \alpha_V)$
 is a short exact sequence of homomorphisms of compatible  Hom-Lie triple systems
$$\begin{CD}%\label{dia:ext}
0@>>> {V} @>i >> \hat{\mathfrak{g}} @>p >> \mathfrak{g} @>>>0\\
@. @V {\alpha_V} VV @V \hat{\alpha}_\mathfrak{g} VV @V \alpha_\mathfrak{g} VV @.\\
0@>>> {V} @>i >> \hat{\mathfrak{g}} @>p >> \mathfrak{g} @>>>0
\end{CD}$$
such that $[u, v, -]^1_{\hat{\mathfrak{g}}}=[u,   v, -]^2_{\hat{\mathfrak{g}}}=0$, for all $u,v\in V$,  i.e., $V$ is an abelian ideal of $\hat{\mathfrak{g}}.$
\end{defn}

A   section  of an abelian extension $(\hat{\mathfrak{g}},[-,-,-]^1_{\hat{\mathfrak{g}}}, [-,-,-]^2_{\hat{\mathfrak{g}}}, \alpha_{\hat{\mathfrak{g}}})$ of $(\mathfrak{g},[-,-,-]^1_\mathfrak{g}, [-,-,-]^2_\mathfrak{g}, \alpha_\mathfrak{g})$ by  $(V,[-,-,-]^1_V, [-,-,-]^2_V, \alpha_V)$ is a linear map $\sigma:\mathfrak{g}\rightarrow \hat{\mathfrak{g}}$ such that $p\circ \sigma=\mathrm{id}_\mathfrak{g}$.

 Note that an abelian extension induces a representation of the compatible  Hom-Lie triple system
 $(\mathfrak{g},[-,-,-]^1_\mathfrak{g}, [-,-,-]^2_\mathfrak{g}, \alpha_\mathfrak{g})$ on $(V,  \alpha_V)$  with the maps
   \begin{align}
\theta^i_V(a,b)(u)=[u,\sigma(a),\sigma(b)]^i_{\hat{\mathfrak{g}}},~~\forall ~a,b\in \mathfrak{g}, u\in V~~ \text{and}~~ i=1,2. \label{5.1}
   \end{align}
Clearly, $D^i_V(a,b)(u)=[\sigma(a),\sigma(b),u]^i_{\hat{\mathfrak{g}}}=\theta^i_V(b,a)(u)-\theta^i_V(a,b)(u)$.
It is easy to see that the above  representation is independent of the choice of the section $\sigma$.

Let $\sigma:\mathfrak{g}\rightarrow \hat{\mathfrak{g}}$  be a section. Define linear maps $\mu_1,\mu_2\in \mathrm{Hom}(\wedge^2\mathfrak{g}\otimes\mathfrak{g},\mathfrak{g})$ by
   \begin{align}
&\mu_i(a,b,c)=[\sigma(a),\sigma(b),\sigma(c)]^i_{\hat{\mathfrak{g}}}-\sigma([a,b,c]^i_\mathfrak{g}),~~\forall ~a,b,c\in \mathfrak{g}  ~~ \text{and}~~ i=1,2.\label{5.2}
   \end{align}
   We transfer the compatible  Hom-Lie triple system structure on $\hat{ \mathfrak{g}}$ to $ \mathfrak{g}\oplus V$ by endowing $ \mathfrak{g}\oplus V$  with two multiplications $[-, -,-]_{\mu_1}$ and  $[-, -,-]_{\mu_2}$  defined by
\begin{align}
[(a,u), (b,v),(c,w)]_{\mu_i}&=([a, b, c]^i_{\mathfrak{g}},\theta^i_V(b,c)u-\theta^i_V(a,c)v+D^i_V(a,b)w+\mu_i(a,b,c)), \label{5.3}
\end{align}
for any $a,b,c\in \mathfrak{g}, u,v,w\in V$  and  $i=1,2.$
Then, it is routine to check that   $(\mathfrak{g}\oplus V,[-, -,-]_{\mu_1},[-, -,-]_{\mu_2},$ $ \alpha_\mathfrak{g}\oplus\alpha_V)$ is a compatible  Hom-Lie triple system  and
$(\mu_1,\mu_2)\in \mathcal{C}_{\mathrm{cHLts}}^{2}( \mathfrak{g},  V)$ is a  2-cocycle in the cohomology of the compatible Hom-Lie triple system  $\mathfrak{g}$ with coefficients in the   representation $(V, \theta^1_V, \theta^2_V, \alpha_V)$.

\begin{defn}
Let $(\hat{\mathfrak{g}}_1,[-,-,-]^1_{\hat{\mathfrak{g}}_1}, [-,-,-]^2_{\hat{\mathfrak{g}}_1}, \alpha_{\hat{\mathfrak{g}}_1})$ and  $(\hat{\mathfrak{g}}_2,[-,-,-]^1_{\hat{\mathfrak{g}}_2}, [-,-,-]^2_{\hat{\mathfrak{g}}_2}, \alpha_{\hat{\mathfrak{g}}_2})$ be two abelian extensions of $(\mathfrak{g},[-,-,-]^1_\mathfrak{g}, [-,-,-]^2_\mathfrak{g}, \alpha_\mathfrak{g})$ by  $(V,[-,-,-]^1_V, [-,-,-]^2_V, \alpha_V)$. They are said to be  equivalent if  there is an isomorphism of  compatible  Hom-Lie triple system  $\zeta: (\hat{\mathfrak{g}}_1,[-,-,-]^1_{\hat{\mathfrak{g}}_1}, [-,-,-]^2_{\hat{\mathfrak{g}}_1}, \alpha_{\hat{\mathfrak{g}}_1})\rightarrow(\hat{\mathfrak{g}}_2,[-,-,-]^1_{\hat{\mathfrak{g}}_2}, [-,-,-]^2_{\hat{\mathfrak{g}}_2}, \alpha_{\hat{\mathfrak{g}}_2})$
such that the following diagram is  commutative:
\begin{align}
\begin{CD}
0@>>> {V} @>i_1 >> \hat{ \mathfrak{g}}_1 @>p_1 >> \mathfrak{g} @>>>0\\
@. @| @V \zeta VV @| @.\\
0@>>> {V} @>i_2 >> \hat{\mathfrak{g}}_2 @>p_2 >>  \mathfrak{g} @>>>0.\label{5.4}
\end{CD}
\end{align}
\end{defn}

Next we are ready to classify abelian extensions of a compatible  Hom-Lie triple systems.

\begin{theorem}
Abelian extensions of a compatible  Hom-Lie triple systems $(\mathfrak{g},[-,-,-]^1_\mathfrak{g},$ $ [-,-,-]^2_\mathfrak{g}, \alpha_\mathfrak{g})$ by  $(V,[-,-,-]^1_V, [-,-,-]^2_V, \alpha_V)$ are classified by the second cohomology group $\mathcal{H}_{\mathrm{cHLts}}^{2}( \mathfrak{g},  V)$ of $(\mathfrak{g},[-,-,-]^1_\mathfrak{g},$ $ [-,-,-]^2_\mathfrak{g}, \alpha_\mathfrak{g})$ with coefficients in the representation $(V, \theta^1_V, \theta^2_V, \alpha_V)$.
\end{theorem}
\begin{proof}
Let $(\hat{\mathfrak{g}},[-,-,-]^1_{\hat{\mathfrak{g}}}, [-,-,-]^2_{\hat{\mathfrak{g}}}, \alpha_{\hat{\mathfrak{g}}})$ be an abelian extension of $(\mathfrak{g},[-,-,-]^1_\mathfrak{g},$ $ [-,-,-]^2_\mathfrak{g}, \alpha_\mathfrak{g})$ by  $(V,[-,-,-]^1_V, [-,-,-]^2_V, \alpha_V)$.  We choose a section $\sigma:  \mathfrak{g}\rightarrow \hat{  \mathfrak{g}}$ to obtain a 2-cocycle $( \mu_1, \mu_2)$.
First, we show that the cohomological class of $( \mu_1, \mu_2)$  is independent of the choice of $\sigma$.
Let $\sigma_1,\sigma_2: \mathfrak{g}\rightarrow \hat{  \mathfrak{g}}$ be two distinct sections providing 2-cocycles $(\mu_1,\mu_2)$ and $(\nu_1,\nu_2)$ respectively. Define
linear map $\xi:  \mathfrak{g}\rightarrow  V$ by $\xi(a)=\sigma_1(a)-\sigma_2(a)$. Then, we get
 $$(\mu_1,\mu_2)-(\nu_1,\nu_2)=(\delta_1\xi,\delta_2\xi)=\delta_c(\xi)\in \mathcal{B}_{\mathrm{cHLts}}^{2}( \mathfrak{g},  V).$$
  So $( \mu_1, \mu_2)$ and $(\nu_1,\nu_2)$ are in the same cohomological class  in $\mathcal{H}_{\mathrm{cHLts}}^{2}( \mathfrak{g},  V)$.

Next, assume that $(\hat{\mathfrak{g}}_1,[-,-,-]^1_{\hat{\mathfrak{g}}_1}, [-,-,-]^2_{\hat{\mathfrak{g}}_1}, \alpha_{\hat{\mathfrak{g}}_1})$ and  $(\hat{\mathfrak{g}}_2,[-,-,-]^1_{\hat{\mathfrak{g}}_2}, [-,-,-]^2_{\hat{\mathfrak{g}}_2}, \alpha_{\hat{\mathfrak{g}}_2})$   are two equivalent abelian extensions   of $(\mathfrak{g},[-,-,-]^1_\mathfrak{g}, [-,-,-]^2_\mathfrak{g}, \alpha_\mathfrak{g})$ by  $(V,[-,-,-]^1_V,$ $ [-,-,-]^2_V, \alpha_V)$ with the associated isomorphism $\zeta: (\hat{\mathfrak{g}}_1,[-,-,-]^1_{\hat{\mathfrak{g}}_1}, [-,-,-]^2_{\hat{\mathfrak{g}}_1},  \alpha_{\hat{\mathfrak{g}}_1})\rightarrow(\hat{\mathfrak{g}}_2,[-,-,-]^1_{\hat{\mathfrak{g}}_2}, [-,-,-]^2_{\hat{\mathfrak{g}}_2}, \alpha_{\hat{\mathfrak{g}}_2})$ such that the diagram in \eqref{5.4} is commutative.  Let $\sigma_1$ be a section of  $(\hat{\mathfrak{g}}_1,[-,-,-]^1_{\hat{\mathfrak{g}}_1}, [-,-,-]^2_{\hat{\mathfrak{g}}_1}, \alpha_{\hat{\mathfrak{g}}_1})$. As $p_2\circ\zeta=p_1$, we  have
$$p_2\circ(\zeta\circ \sigma_1)=p_1\circ \sigma_1= \mathrm{id}_{\mathfrak{g}}.$$
That is, $\zeta\circ \sigma_1$ is a section of $(\hat{\mathfrak{g}}_2,[-,-,-]^1_{\hat{\mathfrak{g}}_2}, [-,-,-]^2_{\hat{\mathfrak{g}}_2}, \alpha_{\hat{\mathfrak{g}}_2})$ . Denote $\sigma_2:=\zeta\circ \sigma_1$. Since $\zeta$ is a isomorphism of  compatible  Hom-Lie triple system such that $\zeta|_V=\mathrm{id}_V$, then, for any $a,b,c\in \mathfrak{g}$, we have
\begin{align*}
\nu_1(a,b,c)&=[\sigma_2(a), \sigma_2(b), \sigma_{2}(c)]^1_{\hat{\mathfrak{g}}_2}-\sigma_2([a,b,c]^1_{\mathfrak{g}_2})\\
&=[\zeta(\sigma_1(a)), \zeta(\sigma_1(b)), \zeta(\sigma_1(c))]^1_{\hat{\mathfrak{g}}_2}-\zeta(\sigma_1([a, b, c]^1_{\mathfrak{g}_2}))\\
&=\zeta\big([\sigma_1(a), \sigma_1(b), \sigma_1(c)]^1_{\hat{\mathfrak{g}}_1}-\sigma_1([a, b, c]^1_{\mathfrak{g}_1})\big)\\
&=\zeta(\mu_1(a,b,c))=\mu_1(a,b,c)
\end{align*}
Similarly, we have $\nu_2=\mu_2.$
Hence, all equivalent abelian extensions give rise to the same element in $\mathcal{H}_{\mathrm{cHLts}}^{2}( \mathfrak{g},  V)$.

Conversely, given two  cohomologous 2-cocycles $(\mu_1,\mu_2)$ and $(\nu_1,\nu_2)$  in $\mathcal{H}_{\mathrm{cHLts}}^{2}( \mathfrak{g},  V)$,
we can construct two abelian extensions $(\mathfrak{g}\oplus V,[-, -,-]_{\mu_1},[-, -,-]_{\mu_2},$ $ \alpha_\mathfrak{g}\oplus\alpha_V)$ and  $(\mathfrak{g}\oplus V,[-, -,-]_{\nu_1},[-, -,-]_{\nu_2},$ $ \alpha_\mathfrak{g}\oplus\alpha_V)$  via Eq.~\eqref{5.3}. Then  there is  a linear map $\xi:  \mathfrak{g}\rightarrow  V$ such that
 $$(\mu_1,\mu_2)-(\nu_1,\nu_2)=(\delta_1\xi, \delta_2\xi)=\delta_c(\xi).$$

 Define linear map $\zeta_\xi:  \mathfrak{g}\oplus  V\rightarrow   \mathfrak{g}\oplus V$ by
$\zeta_\xi(a,u):=(a,\xi(a)+u), ~a\in  \mathfrak{g}, u\in V.$ It is obvious that  $\zeta_\xi$ is an isomorphism of these two abelian extensions such that the diagram in \eqref{5.4} is commutative.
\end{proof}

{{\bf Acknowledgments.}  The paper is  supported by the  Foundation of Science and Technology of Guizhou Province(Grant Nos. [2018]1020,  ZK[2022]031,  QKHZC[2023]372),
the Scientific Research Foundation for Science \& Technology Innovation Talent Team of the Intelligent Computing and Monitoring of Guizhou Province (Grant No. QJJ[2023]063),
 the National Natural Science Foundation of China (Grant No. 12161013).

\end{document}